# Infinity of solutions to initial-boundary value problems for linear constant-coefficient evolution PDEs on semi-infinite intervals

Andreas Chatziafratis [1,3,*] and Spyridon Kamvissis [2,3]

**Abstract.** In this short communication, we announce an algorithmic procedure for constructing non-uniqueness counter-examples of classical solutions to initial-boundary-value problems for a wide class of linear evolution partial differential equations, of any order and with constant coefficients, formulated in a quarter-plane. Our approach relies on analysis of regularity and asymptotic properties, near the boundary of the spatio-temporal domain, of closed-form integral-representation formulae derived via complex-analytic techniques and rigorous implementation of the modern PDE technique known as Fokas unified transform method. In order to elucidate the novel idea and demonstrate the proposed technique in a self-contained fashion, we explicitly present its application to two concrete examples, namely the heat equation and the linear KdV equation with Dirichlet data. New uniqueness theorems for these two models are also presented herein.

**Background.** In order to place this work in context we begin with some introductory remarks. A few decades ago, a technique was pioneered by Fokas [1-3], nowadays broadly known as unified transform method (UTM), for analyzing *integrable* nonlinear PDEs, which was based on the synergy between the *Riemann-Hilbert* problems theory and the *Lax-pairs* formalism. It was later realized [4] that this novel method offered a convenient and advantageous (in contrast to conventional transform methods and their well-known pathologies [4-10], even in the limited cases when the latter are 'applicable') framework for studying linear PDEs posed in semi-bounded domains. In the latter case, the UTM yields solution representations expressed as contour integrals in the spectral complex plane, and is thus regarded as the analogue of the *Fourier* transform which is traditionally involved in the solution of the corresponding problems on the whole space.

More specifically, regarding initial-boundary-value problems (IBVPs), for linear PDEs with *polynomial* dispersion relation, the UTM methodology is an algorithmic procedure for formally deriving formulae (as *candidate* solutions) for such evolution equations posed on a semi-infinite line or a finite interval, assuming existence and certain properties for the sought-after solution which allow the proposed steps of the recipe (as is always the case with transform methods). In particular, the UTM consists of three main steps: (1) Formal construction, assuming validity of all calculations involved, of an integral representation, in the complex plane, for the unknown solution; this representation is not yet effective because it contains certain transforms of unknown boundary values. (2) Analysis of the so-called *global* relation, which is a simple integro-algebraic equation connecting the given initial and boundary data with the unknown boundary values. (3) Use of certain invariance properties of the global relation dictated by the so-called *symmetry* relation combined with algebraic manipulations in order to eliminate the transforms of the unknown boundary values arising in the first step of the algorithm.

This contemporary method has been successful in a broad class of linear-PDE problems (of which the greatest part is *not* amenable to treatment via classical methods) and has therefore been established as an efficient approach, both analytically and numerically, for applications in connection with areas as diverse as boundary value problems, spectral theory, control theory, inverse problems, and so on. For all these areas of application of the UTM, and extensions thereof by many mathematicians, see, for example, [11-28], in chronological order, and more references cited therein. Incidentally, the UTM has also led to the development of a new approach to investigating well-posedness questions for nonlinear PDEs, e.g. [29-32].

However, surprisingly, the literature severely lacks complete and accurate results on existence and uniqueness of classical solutions to IBVPs for linear PDEs. Indicatively, in particular, as mentioned above, the derivation of a (candidate) solution formula for a given IBVP for a PDE proceeds with formal application of Green's theorem to the divergence form of the equation on a semi-infinite strip, which in turn a priori assumes, inter alia, sufficient regularity of the yet unknown solution as well as uniform – in time – decay for the solution and certain derivatives as the spatial variable tends to infinity. Remarkably, the formulae constructed via the UTM were never (until the recent appearance of a new line of rigorous investigations – see below) verified *a posteriori* (which would also guarantee *existence* of solutions). Apparently, this is a common pitfall in the literature dealing with 'exact

[1] Department of Mathematics, National and Kapodistrian University of Athens, Greece

[2] Department of Pure and Applied Mathematics, University of Crete, Greece

[3] Institute of Applied and Computational Mathematics, Foundation for Research and Technology – Hellas, Crete, Greece

*Key words and phrases*: evolution PDEs on semi-bounded domains, nonhomogeneous initial–boundary-value problems, Fokas unified-transform-method, complex-contour integral representations, classical solutions, boundary hehavior, ill-posedness, non-uniqueness.

*2020 Mathematics Subject Classification*: 35A02, 35A09, 35A22, 35A25, 35C05, 35C15, 35G16, 35K57, 35K70, 35Q35, 35Q53, 35Q79, 44A15, 74H25, 74S70, 76B15, 80A19

solutions' via transform methods (traditional or not). Moreover, in general, one cannot know whether the formula arrived at is a *unique* solution, since the derivation process is admittedly formal. Stated differently, even if the derivation were *rigorous* it would still yield *all candidate* explicit solutions; then, rigorous *back-substitution* (or justification) would indeed secure the validity of the proposed closed-form solution.

In view of remedying the aforementioned somewhat-careless reasoning throughout the literature, a new analytical approach was recently introduced in [33], and developed in e.g. [34-43], in alphabetical order, for the rigorous refinement and generalization of the UTM as well as for analytical investigation of a miscellany of qualitative properties of such PDEs in classical settings, e.g. constructive existence and spatiotemporal asymptotics. In the particular case of evolution PDEs with polynomial dispersion relation, the novel approach utilizes, as a starting point, the formulae afforded by suitable implementation of the UTM. One of the salient aspects of these works is the proper *a posteriori* justification of the validity of formally derived UTM representations, including reconstruction –in the sense of limits near the semi-infinite boundaries– of prescribed initial and boundary data – far from a trivial task, in the majority of cases. Notably, these rigorous investigations around the UTM have led to discovery of new instability, blow-up and break-down phenomena; e.g. [38,41,34].

As mentioned already, in this note, we present new examples of *non*-uniqueness for two basic IBVPs (the "simplicity" of these problems make the discovery even more startling). Their justification follows in a straightforward manner from the requisite qualitative analysis published e.g. in [42, 43]. These examples clearly indicate the pattern of a novel 'recipe' for constructing such examples (not to be found in the existing literature) for a large class of equations containing also higher-order evolution PDEs which involve lower-order terms too.

In light of these counter-examples, one way to establish uniqueness is to impose certain data compatibility and solution integrability conditions (which, in a sense, are optimal, as can be illustrated with an explicit calculation, for instance in the case of the heat equation) which in turn permit applicability of an "energy estimate" type of argument, as performed in [37] for the Milne-Barenblatt pseudo-parabolic equation of radiation diffusion and of seepage in porous media. Such propositions are presented for the first time officially in this letter. Further applications to other celebrated models (e.g. for 4th-order equations and coupled systems of PDE) and more related results (pertaining to, e.g., controllability, continuous dependence on data, and how problems with oblique Robin data can effectively be reduced to the Dirichlet-data case) will be reported in subsequent publications.

Although the present line of investigation is based on mostly elementary machinery (thereby also allowing access to a wide readership), yet the derived results appear to be, to the best of our knowledge, both new and surprising. Similar ill-posedness issues for the very specific case of the heat equation have been addressed in historical works of Tychonov [44] and Widder [45], albeit without mention of corresponding consistency conditions on the data at the corner of the domain, i.e. the origin. Namely our findings on the heat equation are deemed complementary to those classical works cited above. However, for higher-order PDE, such as the KdV equation, analogous results are to date apparently missing altogether from the published literature.

**Example 1.** The UTM leads to the following solution:

$$v(x,t) = \frac{i}{\pi}\int_{\gamma}[e^{i\lambda x - \lambda^2 t} - e^{i\lambda x}]\frac{d\lambda}{\lambda}, \text{ defined for } x > 0, \, t > 0,$$

to the IBVP: $v_t = v_{xx}$, $v(0,t) = \lim_{x \to 0^+} v(x,t) = 1$, $v(x,0) = \lim_{t \to 0^+} v(x,t) = 0$. (See e.g. [43].) The contour $\gamma$ is the oriented boundary of the domain $\Theta^- := \{\lambda \in \mathbb{C} : \text{Im}\,\lambda \geq 0, \text{Re}(\lambda^2) \leq 0\}$ (see Fig. 1).

By Cauchy's theorem and Jordan's lemma, this solution can be written also in the following way:

$$v(x,t) = \frac{i}{\pi}\int_{-\infty}^{\infty}[e^{i\lambda x - \lambda^2 t} - e^{i\lambda x}]\frac{d\lambda}{\lambda} = \frac{i}{\pi}\int_{\gamma_0} e^{i\lambda x - \lambda^2 t}\frac{d\lambda}{\lambda},$$

where $\gamma_0 := (\gamma \cap \{|\lambda| \geq \sqrt{2}\}) + [-1+i, 1+i]$ (see Fig.2).

Its derivative with respect to $t$,

$$u(x,t) = \frac{\partial v(x,t)}{\partial t} = -\frac{i}{\pi}\int_{\gamma_0} e^{i\lambda x - \lambda^2 t}\lambda d\lambda = -\frac{i}{\pi}\int_{-\infty}^{\infty} e^{i\lambda x - \lambda^2 t}\lambda d\lambda = \frac{1}{2\sqrt{\pi}}\frac{x}{t\sqrt{t}}e^{-x^2/4t},$$

(which – upto a constant – is the $x$-derivative of Gauss's kernel), solves the IBVP:

$$u_t = u_{xx} \; (x > 0, t > 0), \; u(0,t) = \lim_{x \to 0^+} u(x,t) = 0 \; (t > 0), \; u(x,0) = \lim_{t \to 0^+} u(x,t) = 0 \; (x > 0), \tag{1}$$

and $u(x,t) \not\equiv 0$.

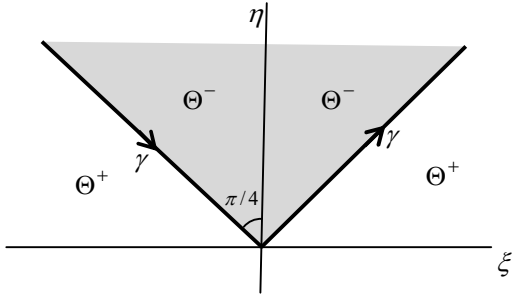 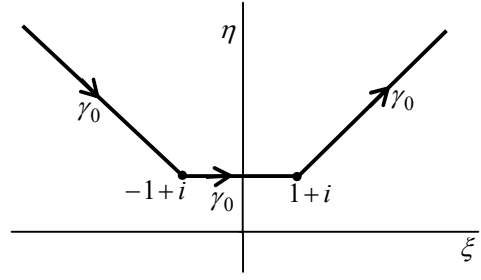

**Fig.1** *The contour $\gamma = \partial\Theta^-$ in the $\lambda$-plane ($\lambda = \xi + i\eta$)*   **Fig.2** *The contour $\gamma_0$*

Thus, the problem (1) does not have unique solution. As a matter of fact, each of the functions

$$u_n(x,t) := \frac{\partial^n v(x,t)}{\partial t^n} = -\frac{i}{\pi}\int_{-\infty}^{\infty}(-\lambda^2)^{n-1}e^{i\lambda x - \lambda^2 t}\lambda d\lambda = \frac{x}{2\sqrt{\pi}}\frac{\partial^{n-1}}{\partial t^{n-1}}\left(\frac{1}{t\sqrt{t}}e^{-x^2/4t}\right), \quad n = 1,2,3,\ldots,$$

solves (1).

**Example 2.** The UTM yields [42] the following solution formula:

$$v(x,t) = -\frac{3}{2\pi i}\int_{\text{Im}\,\lambda=\varepsilon} e^{i\lambda x + i\lambda^3 t}\frac{d\lambda}{\lambda}, \text{ defined for } x > 0, \ t > 0,$$

to the IBVP: $v_t = -v_{xxx}$, $v(0,t) = \lim_{x\to 0^+} v(x,t) = 1$, $v(x,0) = \lim_{t\to 0^+} v(x,t) = 0$. (See [42] for the rigorous justification and full boundary behavior analysis.) The above integral does not depend on $\varepsilon$, where $\varepsilon$ is any fixed positive number. Its derivative with respect to $t$,

$$u(x,t) = \frac{\partial v(x,t)}{\partial t} = -\frac{9}{2\pi}\int_{\text{Im}\,\lambda=\varepsilon}e^{i\lambda x + i\lambda^3 t}\lambda^2 d\lambda = \frac{9}{2\pi}\frac{\partial^2}{\partial x^2}\left(\int_{\text{Im}\,\lambda=\varepsilon}e^{i\lambda x + i\lambda^3 t}d\lambda\right),$$

solves the IBVP:

$$u_t = -u_{xxx} \ (x > 0, t > 0), \ u(0,t) = \lim_{x\to 0^+} u(x,t) = 0 \ (t > 0), \ u(x,0) = \lim_{t\to 0^+} u(x,t) = 0 \ (x > 0). \quad (2)$$

We claim also that $u(x,t) \not\equiv 0$. Indeed, if $u(x,t) \equiv 0$ then $\frac{\partial v(x,t)}{\partial t} \equiv 0$, and this would imply $v(x,t) \equiv v(x,0) \equiv 0$, which would contradict $v(0,t) = 1$.

Also, each of the functions

$$u_n(x,t) := \frac{\partial^n v(x,t)}{\partial t^n} = -\frac{9}{2\pi}\int_{\text{Im}\,\lambda=\varepsilon}(i\lambda^3)^{n-1}e^{i\lambda x + i\lambda^3 t}\lambda^2 d\lambda, \quad n = 1,2,3,\ldots,$$

solves problem (2).

(We refer again to Chatziafratis-Kamvissis-Stratis [42] for properties of the derivatives of the linear KdV solution.)

Now, we consider the following cases of IBVPs with the purpose of stating our precise uniqueness propositions.

**Problem 1.** *Solve*

$$\begin{cases} \dfrac{\partial U}{\partial t} + \dfrac{\partial^3 U}{\partial x^3} = f, \ (x,t) \in Q := \mathbb{R}^+ \times \mathbb{R}^+ \\ \lim_{t\to 0^+} U(x,t) = u_0(x), \ x \in \mathbb{R}^+ \\ \lim_{x\to 0^+} U(x,t) = g_0(t), \ t \in \mathbb{R}^+, \end{cases} \quad (3)$$

*for $U = U(x,t)$.*

**Remark.** *We assume that*

$$u_0 = u_0(x) \in \mathcal{S}([0,\infty)), \ g_0 = g_0(t) \in C^\infty([0,\infty)) \ \text{and} \ f = f(x,t) \in \mathcal{S}(\overline{Q}). \tag{4}$$

Setting

$$\hat{u}_0(\lambda) = \int_{y=0}^{\infty} e^{-i\lambda y} u_0(y) dy \ \text{and} \ \hat{f}(\lambda,t) = \int_{y=0}^{\infty} e^{-i\lambda y} f(y,t) dy, \ \text{both defined for} \ \lambda \in \mathbb{C} \ \text{with} \ \text{Im}\,\lambda \le 0,$$

$$\tilde{g}_0(\omega(\lambda),t) = \int_{\tau=0}^{t} e^{\omega(\lambda)\tau} g_0(\tau) d\tau, \ \text{where} \ \omega(\lambda) = -i\lambda^3 \ (\lambda \in \mathbb{C}),$$

and

$$\tilde{\hat{f}}(\lambda,\omega(\lambda),t) = \int_{\tau=0}^{t} e^{\omega(\lambda)\tau} \hat{f}(\lambda,\tau) d\tau \ (\lambda \in \mathbb{C} \ \text{with} \ \text{Im}\,\lambda \le 0),$$

we define

$$(\mathbb{J}_\circ^+ u_0)(x,t) = \int_{\lambda=-\infty}^{\infty} e^{i\lambda x - \omega(\lambda)t} \hat{u}_0(\lambda) d\lambda, \quad (\mathbb{J}_\circ^- u_0)(x,t) = \int_{\lambda \in \Gamma} e^{i\lambda x - \omega(\lambda)t} [\alpha \hat{u}_0(\alpha\lambda) + \alpha^2 \hat{u}_0(\alpha^2\lambda)] d\lambda,$$

$$(\mathbb{J}_1 g_0)(x,t) = \int_{\lambda \in \Gamma} e^{i\lambda x - \omega(\lambda)t} 3\lambda^2 \tilde{g}_0(\omega(\lambda),t) d\lambda, \quad (\mathbb{J}_2^+ f)(x,t) = \int_{\lambda=-\infty}^{\infty} e^{i\lambda x - \omega(\lambda)t} \tilde{\hat{f}}(\lambda,\omega(\lambda),t) d\lambda,$$

and $\quad (\mathbb{J}_2^- f)(x,t) = \int_{\lambda \in \Gamma} e^{i\lambda x - \omega(\lambda)t} [\alpha \tilde{\hat{f}}(\alpha\lambda,\omega(\lambda),t) + \alpha^2 \tilde{\hat{f}}(\alpha^2\lambda,\omega(\lambda),t)] d\lambda,$

where $\alpha = e^{2\pi i/3}$ and $\Gamma = \partial\Omega^-$ with $\Omega^- = \{\lambda \in \mathbb{C} : \text{Im}\,\lambda \ge 0 \ \text{and} \ \text{Re}\,\omega(\lambda) \le 0\}$ (see Fig. 3).

With this notation, the unified transform method (UTM) leads [42] to the solution of the above problem (3) in the following form: For $x > 0$ and $t > 0$,

$$2\pi U(x,t) = (\mathbb{J}_\circ^+ u_0)(x,t) + (\mathbb{J}_\circ^- u_0)(x,t) - (\mathbb{J}_1 g_0)(x,t) + (\mathbb{J}_2^+ f)(x,t) + (\mathbb{J}_2^- f)(x,t). \tag{5}$$

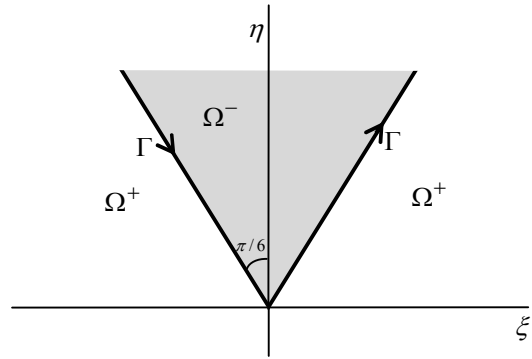

**Fig. 3** *The contour $\Gamma$ is the boundary of $\Omega^-$*

**Theorem 1.** *The function $U(x,t)$, defined by (5), is the unique solution of (3), in the following sense: If, in addition to (4), we assume that $u_0(0) = g_0(0)$ and $g_0'(0) = -u_0'''(0) + f(0,0)$, and*

$$V(x,t) \ \text{is a} \ C^3 \ \text{function in} \ \overline{Q} - \{(0,0)\} \ \text{and solves (3)},$$

$$\lim_{x\to\infty} V(x,t) = \lim_{x\to\infty} V_x(x,t) = 0 \ \text{and} \ \sup_{x\ge 1} |V_{xx}(x,t)| < \infty \ (\forall t > 0),$$

*and, for every $T > 0$, the functions*

$$|V(x,t)|^2, \ |V_t(x,t)|^2 \ \text{are, uniformly for} \ 0 < t \le T, \ \text{integrable with respect to} \ x \in [0,\infty),$$

*i.e., there exists a positive function $B_T(x)$ such that $\int_0^\infty B_T(x) dx < +\infty$ and, for $0 < t \le T$,*

$$|V(x,t)|^2 \le B_T(x) \ \text{and} \ |V_t(x,t)|^2 \le B_T(x) \ (\forall x > 0),$$

*then $V \equiv U$.*

**Problem 2.** *Assuming (1.4), solve*

$$\begin{cases} \dfrac{\partial U}{\partial t} - \dfrac{\partial^2 U}{\partial x^2} = f, \ (x,t) \in Q := \mathbb{R}^+ \times \mathbb{R}^+ \\ \lim_{t \to 0^+} U(x,t) = u_0(x), \ x \in \mathbb{R}^+ \\ \lim_{x \to 0^+} U(x,t) = g_0(t), \ t \in \mathbb{R}^+, \end{cases} \quad (6)$$

*for* $U = U(x,t)$.

For $x > 0$ and $t > 0$, the UTM yields [43] the following solution representation:

$$2\pi U(x,t) = \int_{-\infty}^{\infty} e^{i\lambda x - \lambda^2 t} \hat{u}_0(\lambda) d\lambda - \int_{\gamma} e^{i\lambda x - \lambda^2 t} \hat{u}_0(-\lambda) d\lambda - 2i \int_{\gamma} e^{i\lambda x - \lambda^2 t} \lambda \tilde{g}_0(\omega(\lambda), t) d\lambda$$

$$+ \int_{-\infty}^{\infty} e^{i\lambda x - \lambda^2 t} \tilde{\hat{f}}(\lambda, \omega(\lambda), t) d\lambda - \int_{\gamma} e^{i\lambda x - \lambda^2 t} \tilde{\hat{f}}(-\lambda, \omega(\lambda), t) d\lambda, \quad (7)$$

where, with $\omega(\lambda) = \lambda^2$,

$$\tilde{g}_0(\omega(\lambda), t) = \int_{\tau=0}^{t} e^{\omega(\lambda)\tau} g_0(\tau) d\tau \ (\lambda \in \mathbb{C}) \ \text{and} \ \tilde{\hat{f}}(\lambda, \omega(\lambda), t) = \int_{\tau=0}^{t} e^{\omega(\lambda)\tau} \hat{f}(\lambda, \tau) d\tau \ (\lambda \in \mathbb{C} \ \text{with} \ \text{Im} \lambda \leq 0),$$

the contour $\gamma$ is the oriented boundary of the domain $\{\lambda \in \mathbb{C} : \text{Im} \lambda \geq 0 \ and \ \text{Re}(\lambda^2) \leq 0\}$ (see Fig.1) and $\hat{u}_0(\lambda)$ and $\hat{f}(\lambda, \tau)$ are as defined above.

**Theorem 2.** *The function* $U(x,t)$, *defined by (7), is the unique solution of (6), in the following sense: If, in addition to (4), we assume that*

$$u_0(0) = g_0(0) \ and \ u_0''(0) + f(0,0) = g_0'(0), \quad (8)$$

*and*

$$V(x,t) \ is \ a \ C^2 \ function \ in \ \overline{Q} - \{(0,0)\} \ and \ solves \ (6),$$

$$\lim_{x \to \infty} V(x,t) = 0 \ and \ \sup_{x \geq 1} |V_x(x,t)| < \infty \ (\forall t > 0),$$

*and, for every* $T > 0$,

$$|V(x,t)|^2, \ |V_t(x,t)|^2 \ are, \ uniformly \ for \ 0 < t \leq T, \ integrable \ with \ respect \ to \ x \in [0, \infty),$$

*i.e., there exists a positive function* $B_T(x)$ *such that* $\int_0^{\infty} B_T(x) dx < +\infty$ *and, for* $0 < t \leq T$,

$$|V(x,t)|^2 \leq B_T(x) \ and \ |V_t(x,t)|^2 \leq B_T(x) \ (x > 0), \quad (9)$$

*then* $V \equiv U$.

The proofs of the Theorems 1 and 2, above, follow easily from the analogous Theorem 2.4 in [37].

**Comment** The compatibility conditions at the origin (8), imposed in Theorem 2, imply [43] the existence of the limits

$$\lim_{\overline{Q} \ni (x,t) \to (0,0)} U(x,t), \ \lim_{\overline{Q} \ni (x,t) \to (0,0)} U_x(x,t), \ \lim_{\overline{Q} \ni (x,t) \to (0,0)} U_{xx}(x,t) \ and \ \lim_{\overline{Q} \ni (x,t) \to (0,0)} U_{xxx}(x,t),$$

and this, in turn, implies that $|U(x,t)|^2$, $|U_t(x,t)|^2$ are, uniformly for $0 < t \leq T$, integrable with respect to $x \in [0, \infty)$. This is crucial for the proof of the uniqueness assertion.

A specific case of this situation is the function

$$u(x,t) = \frac{1}{2\sqrt{\pi}} \frac{x}{t\sqrt{t}} e^{-x^2/4t}, \ x > 0,$$

of Example 1 which does not satisfy the first of conditions (9), i.e., it fails to be $L^2$ – integrable (with respect to $x$) **uniformly** for $t$ close to $0$. Indeed, for fixed $t > 0$,

$$\int_0^\infty \left(\frac{x}{t\sqrt{t}}e^{-x^2/4t}\right)^2 dx = \frac{1}{t^3}\int_0^\infty x^2 e^{-x^2/2t}dx = -\frac{1}{t^2}\int_0^\infty x\, d(e^{-x^2/2t})$$

$$= \frac{1}{t^2}\int_0^\infty e^{-x^2/2t}dx = \frac{\sqrt{2}}{t^{3/2}}\int_0^\infty e^{-x^2/2t}d(x/\sqrt{2t}) = \sqrt{\frac{\pi}{2}}\frac{1}{t^{3/2}} \quad (t > 0).$$

Let us recall that the above function $u(x,t)$ is defined to be the derivative $v_t(x,t)$ of the UTM solution $v(x,t)$ of the IBVP: $v_t = v_{xx}$, $v(0,t) = 1$, $v(x,0) = 0$. The data in this problem, namely "$v(0,t) = 1$" and "$v(x,0) = 0$" do not satisfy (8), and this explains the fact that the function $v_t(x,t) = u(x,t)$ is not $L^2$ – integrable, uniformly for $t$ close to $0$, as the previous explicit calculation shows.

An analogous remark, concerning the compatibility conditions at the origin stated in Theorem 1 and the existence of limits of the UTM solution of (3) and certain of its derivatives as $(x,t) \to (0,0)$, can be made also in this case. (See [42].)

*Acknowledgments.* Andreas Chatziafratis gratefully acknowledges partial funding, at different stages of this project, from: the Hellenic State Scholarships Foundation (IKY), a postdoctoral research grant awarded by the Hellenic Foundation for Research and Innovation, and the European Research Council (ERC grant 101078061 SINGinGR). A.C. also wishes to express his thankfulness to Professors: E. C. Aifantis, P. Biler, J. L. Bona, B. Deconinck, A. S. Fokas, G. Fournodavlos, L. Grafakos, T. Hatziafratis, A. A. Himonas, G. Karali, P. G. Kevrekidis, A. Miranville, T. Ozawa, D. A. Smith, I. G. Stratis, and S. F. Tian, for useful discussions and for providing encouragement, inspiration and academic support.